\documentclass[sn-mathphys-num]{sn-jnl}% Math and Physical Sciences Numbered Reference Style 
\usepackage{graphicx}%
\usepackage{multirow}%
\usepackage{amsmath,amssymb,amsfonts}%
\usepackage{amsthm}%
\usepackage{mathrsfs}%
\usepackage[title]{appendix}%
\usepackage{xcolor}%
\usepackage{textcomp}%
\usepackage{manyfoot}%
\usepackage{booktabs}%
\usepackage{algorithm}%
\usepackage{algorithmicx}%
\usepackage{algpseudocode}%
\usepackage{listings}%
\usepackage{geometry}%
\usepackage{pgfplots}%
\theoremstyle{thmstyleone}%
\newtheorem{thm}{Theorem}%  meant for continuous numbers
\newtheorem{prop}[thm]{Proposition}% 
\newtheorem{lem}[thm]{Lemma}

\theoremstyle{thmstyletwo}%
\newtheorem{rem}{Remark}%

\theoremstyle{thmstylethree}%

\usepackage{mathtools}
\mathtoolsset{showonlyrefs}
\usepackage{nicematrix}

\newcommand{\BB}{{\bf B}}   
\newcommand{\emptytree}{\varnothing}
\newcommand{\Bg}{\text{\sf B}}

\newcommand{\Ba}{\text{\sf B}^*}

\newcommand{\f}{\text{\hspace{0.3pt}\bf f\hspace{0.4pt}}}

\newcommand{\boldk}{\text{\hspace{0.3pt}\bf k\hspace{0.4pt}}}
\newcommand{\boldK}{\text{\hspace{0.3pt}\bf K\hspace{0.4pt}}}
\newcommand{\D}{\text{\rm D}}
\newcommand{\F}{\text{\hspace{0.3pt}\bf F\hspace{0.4pt}}}

\newcommand{\tree}{\protect\scalebox{1.05}{%
\text{\hspace{0.3pt}\rm  t\hspace{0.3pt}}}}

\renewcommand{\tree}{%
\text{\hspace{0.3pt}\rm  t\hspace{0.3pt}}}

\newcommand{\Tree}{\protect\scalebox{1.0}{%
\text{\hspace{0.5pt}\rm  T\hspace{0.5pt}}}}

\renewcommand{\Tree}{%
\text{\hspace{0.3pt}\rm  T\hspace{0.3pt}}}

\newcommand{\forest}{%
\text{\hspace{0.3pt}\rm  f\hspace{0.3pt}}}

\newcommand{\Forest}{\text{\rm F}}

\newcommand{\forestS}{\text{\scriptsize\sf\bfseries F}}
\newcommand{\ForestS}{\text{\small\sf\bfseries F}}

\usepackage{tikz}
\usepackage{pgfplotstable}
\usetikzlibrary{arrows,positioning,plotmarks,external,patterns,angles,
decorations.pathmorphing,backgrounds,fit,shapes,calc}

\newcounter{ncount}
\newcounter{scount}

\def\metrics#1#2#3{%
\def\treesize{#1}\def\thick{#2}\def\rad{#3}}

\def\vertex#1#2{%
\addtocounter{ncount}{1}
\setcounter{scount}{\thencount}
\addtocounter{scount}{-#1}
\node (\thencount) at ($(\thescount) + ( #2,1)$){};
\draw (\thencount)--(\thescount);
\fill (\thencount) circle;
}

\def\treenv#1{%
\begin{tikzpicture}[x=\treesize mm,y=\treesize mm,radius=\rad pt,line width=\thick pt,inner sep=0,baseline=-0.02cm]
\node (0) at (0,0) {}; \fill (0) circle;
\setcounter{ncount}{0}
#1
\end{tikzpicture}}

\renewcommand{\frame}[3][white]{\draw[#1,thin]  (#2) rectangle (#3);}
\metrics{2.5}{1}{2}

\begin{document}

\title[On the B-series composition theorem]{On the B-series composition theorem}


\author[1]{\fnm{John} C. \sur{Butcher}}\email{butcher@math.auckland.ac.nz}

\author[2]{\fnm{Taketomo} \sur{Mitsui}}\email{tom.mitsui@nagoya-u.jp}

\author*[3]{\fnm{Yuto} \sur{Miyatake}}\email{yuto.miyatake.cmc@osaka-u.ac.jp}

\author[4]{\fnm{Shun} \sur{Sato}}\email{shun@mist.i.u-tokyo.ac.jp}

\affil*[1]{\orgdiv{Department of Mathematics}, \orgname{University of Auckland}, \orgaddress{\street{38 Princes Street}, \city{Auckland}, \postcode{1010}, \country{New Zealand}}}

\affil[2]{\orgdiv{Department of Mathematical Informatics}, \orgname{Graduate School of Informatics, Nagoya University}, \orgaddress{\street{Furo-cho, Chikusa-ku}, \city{Nagoya}, \postcode{464-8601}, \state{Aichi}, \country{Japan}}}

\affil[3]{\orgdiv{Cybermedia Center}, \orgname{Osaka University}, \orgaddress{\street{1-32 Machikaneyama}, \city{Toyonaka}, \postcode{560-0043}, \state{Osaka}, \country{Japan}}}

\affil[4]{\orgdiv{Department of Mathematical Informatics}, \orgname{Graduate School of Information Science and Technology, The University of Tokyo}, \orgaddress{\street{Hongo}, \city{Bunkyo-ku}, \postcode{113-8656}, \state{Tokyo}, \country{Japan}}}

\abstract{The B-series composition theorem has been an important topic in numerical analysis of ordinary differential equations for the past-half century.
    Traditional proofs of this theorem rely on labelled trees, whereas recent developments in B-series analysis favour the use of unlabelled trees. 
    In this paper, we present a new proof of the B-series composition theorem that does not depend on labelled trees. 
    A key challenge in this approach is accurately counting combinations related to ``pruning.'' 
    This challenge is overcome by introducing the concept of ``assignment.''}




\maketitle

\section{Introduction}

Given an $n$-dimensional initial value problem
\begin{equation*}
    \frac{\mathrm{d}}{\mathrm d x}y (x)= f(y(x)), \quad y(x_0) = y_0,
\end{equation*}
where $f\in C^\infty(\mathbb{R}^n,\mathbb{R}^n)$,
B-series play an indispensable role
in the analysis of numerical methods, such as Runge--Kutta methods.
The importance of B-series was recognised in the early 1970s and led to major developments in algebraic analysis.
The composition rule was initially proposed by Butcher in~\cite{bu72}, and a more direct proof was provided by Hairer and Wanner~\cite{hw74} (see also Theorem II.12.6 in~\cite{hnw93}). 
A proof with the more modern notation is also given in~\cite{hlw06}.
These proofs make extensive use of labelled trees.

A recent monograph by the first author in 2021 \cite{bu21} attempted to establish the composition theorem without the use of labelled trees. 
After extensive discussion while translating the book into Japanese, a new proof was found based on a correction to the proof in~\cite{bu21}.
This translation has been published~\cite{bu24}, though the new proof has remained confined to the Japanese language until now. 
This paper seeks to present the new proof in English.

\subsection{B-series composition theorem and a key lemma}

Note that we are using some notation loosely here; precise definitions will be provided at the end of this section.

Let 
\begin{equation}
\label{eq:tree}
\Tree =\bigg\{ \raisebox{-6pt}{\treenv{},\,
    \treenv{\vertex10}, \, 
    \treenv{\vertex1{-0.6}\vertex2{0.6}}, \,
    \treenv{\vertex10\vertex10}, \,
    \treenv{\vertex1{-0.7}\vertex20\vertex3{0.7}}, \,
    \treenv{\vertex1{-0.6}\vertex2{0.6}\vertex10}, \,
    \treenv{\vertex10\vertex1{-0.6}\vertex2{0.6}}, \,
    \treenv{\vertex10\vertex10\vertex10},\, \dots} \bigg\}
\end{equation}
be the set of rooted trees and $\Ba$ as the set of all mappings from $\emptytree \cup \Tree$ to $\mathbb{R}$.
For $a\in \Ba$, 
the B-series $(\BB_h y_0) a$ is a formal series defined by
\begin{equation}\label{eq:3.4.b}
(\BB_hy_0)a = a(\emptytree) y_0 + \sum_{\tree\in \Tree} \frac{h^{|t|}a(\tree)}{\sigma(\tree)}\F(\tree),
\end{equation}
where the elementary differential $\F(\tree)$ is implicitly evaluated at $y_0$.
In particular, when $a(\emptytree)=1$, this formal series not only encompasses the Taylor series expansion for the exact solution $y(x_0 + h)$ but also includes approximate solutions obtained by applying important classes of numerical methods, such as Runge--Kutta methods.
By using the affine subspace $\Bg = \{ a \in \Ba \mid a(\emptytree) = 1\}$, the B-series composition theorem, which is the focus of this paper, can be summarized as follows.

% Note that the notation used in the paper is summarized at the end of this section.

% This paper focuses on a key theorem regarding the composition of two B-series:

\begin{thm}
\label{theorem1}
Let $a\in\Bg$ and $b\in\Ba$. Then,
$\big( \BB_h (\BB_hy_0)a \big) b$
is also a B-series $(\BB_hy_0)(ab)$, and the map $ab \in \Ba$ is given as
\begin{align}
(ab)(\emptytree) &= b(\emptytree), \nonumber \\
(ab)(\tree) &= b(\emptytree) a(\tree) +
 \sum_{\tree'\le\tree}b(\tree') a(\tree\setminus\tree'),\quad \tree\in \Tree. \label{eq:3.9.hh}
\end{align}
\end{thm}

Once the following lemma has been proved,
the proof of Theorem~\ref{theorem1} is straightforward.

\begin{lem}
\label{lemma1}
Let $a\in\Bg$.
For $\tree'\in\Tree$, it follows that
\begin{equation}\label{eq:3.9.h}
\frac{h^{|\tree'|}}{\sigma(\tree')}\F(\tree')\big((\BB_h a)y_0\big) 
=  \sum_{\tree\ge\tree'}    \frac{h^{|\tree|}}{\sigma(\tree)}a(\tree\setminus\tree')\F(\tree) (y_0).
\end{equation}
\end{lem}

The main contribution of the paper is to provide a new proof for Lemma~\ref{lemma1}.
The major challenge in avoiding labelled trees arises from the complexity of counting trees.
To clarify this issue, we introduce the concept of an \emph{assignment} of two trees, which plays a pivotal role in the counting process.

\subsection{Notation}

A graph $ G = (V,E) $ is a pair of the sets $V$ and $E $, where $ V $ is the finite set of vertices and $E \subseteq V \times V $ is the set of edges. 
A (rooted) tree $ \tree = (V,E,r) $ is a triplet such that $ (V,E)$ is a connected graph with at least one vertex and with no loops, and $ r \in V $. 
We identify two trees $ \tree = (V,E,r) $ and $ \tree' = (V',E',r')$ if there exists a bijective map $ \pi \colon V \to V' $ such that $ (v_1, v_2) \in E \iff (\pi(v_1),\pi(v_2)) \in E' $. 
% Based on this identification, $ \Tree $ denotes the set of trees, and the first few elements of the set can be expressed as follows: 
% \begin{align*}
%     &\treenv{},&
%     &\treenv{\vertex10},&
%     &\treenv{\vertex1{-0.6}\vertex2{0.6}},&
%     &\treenv{\vertex10\vertex10},&
%     &\treenv{\vertex1{-0.7}\vertex20\vertex3{0.7}},&
%     &\treenv{\vertex1{-0.6}\vertex2{0.6}\vertex10},&
%     &\treenv{\vertex10\vertex1{-0.6}\vertex2{0.6}},&
%     &\treenv{\vertex10\vertex10\vertex10},
% \end{align*}
% where the lowest vertex of each tree corresponds to the root. 
Recalling the introduction of trees given by \eqref{eq:tree}, we denote the order of $\tree$, that is the number of vertices, as $|\tree|$. 
The unique tree with order $1$ is denoted by $\tau$.
% A tree with exactly one vertex (and no edge) is denoted by $\tau$. 
% We often denote by $|\tree|$ the number of vertices of a tree $\tree$.
A new tree obtained by connecting the roots of $\tree_1,\tree_2,\dots,\tree_m\in \Tree$
with a new root
is denoted by 
\begin{equation*}
\tree = [\tree_1\tree_2\dots\tree_m].
\end{equation*}
If some of $\tree_1,\tree_2,\dots,\tree_m$ are equal to each other, we write 
\begin{equation*}
\tree = [\tree_1^{k_1}\tree_2^{k_2}\dots\tree_m^{k_m}],
\end{equation*}
for example $\Big[\treenv{\vertex10}, \treenv{\vertex1{-0.6}\vertex2{0.6}}, \treenv{\vertex1{-0.6}\vertex2{0.6}}\Big] = \Big[\treenv{\vertex10}, {\treenv{\vertex1{-0.6}\vertex2{0.6}}\, }^2\Big]$.
We introduce $a \D:\emptytree \cup \Tree \to \mathbb{R}$ by generalizing $a \in \Bg $:
\begin{align*}
    & (a\D) (\emptytree) = 0, \quad (a\D) (\tau) = 1, \\
    & (a\D) (\tree) = \prod_{i=1}^m a(\tree_i), \quad \tree = [\tree_1,\dots, \tree_m].
\end{align*}

The function $\sigma:\Tree \to \mathbb{R}$, called the symmetry, is defined recursively by
\begin{align}
\sigma(\tau) &= 1,\nonumber\\
\sigma(\tree)&=
\prod_{i=1}^m k_i! \sigma(\tree_i)^{k_i}, \qquad \tree = [\tree_1^{k_1} \tree_2^{k_2}\cdots \tree_m^{k_m}] .\label{eq:2.5.b}
\end{align}
The elementary differential associated with a tree $\tree=[\tree_1 \tree_2 \cdots \tree_n]$ is defined by 
\begin{align}
\F(\tau)&=\f, \\
\F(\tree)&=\f^{(n)}\F(\tree_1)\F(\tree_2)\cdots \F(\tree_n), \quad \tree = [\tree_1,\dots, \tree_m].
\end{align}
Here, bold notation indicates the omission of $(y_0)$, i.e. $\f = f(y_0)$ and $\F (\tree) = \F (\tree) (y_0)$.
If another vector is taken as an argument, it will be written explicitly.
Additionally, $\f^{(n)}$ denotes the Fréchet derivative as a multi-linear operator evaluated at $y_0$.
For later use, we slightly generalize the definition of elementary differentials. 
For the tree $ \tree $, 
we introduce the concept of a  ``stump'' $ \tree \ast^n $, which is a tree-like object with $n$ edges connected to the root
with the other ends of these edges left blank.
We regard the stump $ \tree \ast^n $ as a mapping $ \tree \ast^n \colon \Tree^n \to \Tree $ satisfying
\[ \tree \ast^n \tree'_1 \tree'_2 \cdots \tree'_n = [ \tree_1 \tree_2 \cdots \tree_m \tree'_1 \tree'_2 \cdots \tree'_n ] ,\]
where $ \tree = [ \tree_1 \tree_2 \cdots \tree_m ] $ (see \cite{bu21} for details on the stumps).
For a tree $ \tree = [\tree_1 \tree_2 \cdots \tree_m ] $, we define
\[ \F (\tree \ast^n ) = \f^{(m+n)} \F (\tree_1 ) \F (\tree_2) \cdots \F (\tree_m). \]

A forest is a juxtaposition of trees. 
The set of forests $\Forest$ can be defined 
as the commutative semigroup with identity $1$, generated by $\Tree$ so that it satisfies the following relations:
\begin{align*}
1 &\in \Forest,\\
\tree\;\forest = \forest \;\tree\ &\in \Forest, \qquad \tree\in\Tree, \quad\forest\in \Forest.
\end{align*}
We denote by $ |\forest|$ 
the order of $\forest$ equal to 
the total of the orders of the trees comprising $\forest$, so that $|[\forest]| = 1 + |\forest|$.
The forest space $ \ForestS$ is a space of formal linear combination of forests, i.e. an element $ \forestS $ of $\ForestS $ can be written as $ \forestS = \sum_{ \forest \in \Forest } c(\forest) \forest$, where $ c \colon \Forest \to \mathbb{R} $. 
The forest space $ \ForestS$ is a unitary commutative ring with the obvious addition and multiplication. 
Then, for $ a \in \Bg $, there exists a unique homomorphic extension to the forest space, which is again denoted by $ a \colon \ForestS \to \mathbb{R} $. 
For example,
\begin{equation*}
    a \Big(\treenv{\vertex1{-0.6}\vertex2{0.6}} + 2 \ \treenv{} \, \treenv{\vertex10} + 3 \ \treenv{}^3 \Big) = 
    a \Big(\treenv{\vertex1{-0.6}\vertex2{0.6}}\Big) + 2 a( \treenv{} ) a\big( \treenv{\vertex10}\big) + 3 a( \treenv{})^3.
\end{equation*}

% \memo{Definition of prunings}

A tree $ \tree' = (V', E', r' ) $ is called a subtree of another tree $ \tree = (V,E,r) $ if there exists an injective mapping $ \pi \colon V' \to V $ such that $ \pi (r') = r $ and $ (\pi(v'_1), \pi (v'_2)) \in E $ for all $ (v'_1, v'_2 ) \in E' $. 
% If $ \tree' $ is a subtree of $ \tree $, $ \tree $ is said to be a supertree of $ \tree' $. 
We denote $ \tree' $ being a subtree of $ \tree $ as $ \tree' \le \tree $ or $ \tree \ge \tree' $.
% It is denoted by $ \tree' \le \tree $ or $ \tree \ge \tree' $ that $ \tree' $ is a subtree of $ \tree $. 

For the pair of trees $ (\tree,\tree') $ such that $ \tree \ge \tree' $, the pruning $ \tree \setminus \tree' $ is defined as a formal sum of the possible ways in which $\tree$ can be pruned to form $\tree'$. 
To define the pruning formally, we introduce the set 
\[ \Pi (\tree',\tree) \coloneqq \{ \pi \colon V' \to V \mid \pi \text{ is injective, } \pi (r') = r , \ (v'_1,v'_2) \in E' \Rightarrow (\pi(v'_1),\pi(v'_2)) \in E \}. \]
Then, the pruning is given by 
\[ \tree \setminus \tree' \coloneqq \sum_{ \pi \in \Pi (\tree',\tree) } G ( V \setminus \pi (V') ), \]
where $ G ( V \setminus \pi (V') ) $ denotes the induced subgraph of the vertex set $ V \setminus \pi (V') $,
i.e. $ G ( V \setminus \pi (V') ) $ is the graph with the vertex set $ V \setminus \pi (V') $ and the edge set $ \{ (v_1,v_2) \mid v_1, v_2 \in V \setminus \pi (V'), \ (v_1,v_2) \in E \} $. 
The following are examples of pruning:
\begin{align*}
\treenv{\vertex1{-0.6}\vertex2{0.6}\vertex10} \setminus \treenv{\vertex10\vertex10} &= \treenv{} \ , &
\treenv{\vertex1{-0.6}\vertex2{0.6}} \setminus \treenv{\vertex10} &= 2 \ \treenv{} \ ,&
\treenv{\vertex1{-0.8}\vertex2{0}\vertex3{0.8}\vertex10\vertex30} \setminus \treenv{\vertex1{-0.6}\vertex2{0.6}\vertex10} &= 2 \ \treenv{\vertex10} + 2 \ \treenv{}^2,&
\treenv{\vertex1{-0.9}\vertex20\vertex10\vertex4{0.9}\vertex1{-0.3}\vertex2{0.3}} \setminus \treenv{\vertex1{-0.6}\vertex2{0.6}\vertex10} &= \treenv{\vertex1{-0.6}\vertex2{0.6}} + 2 \ \treenv{} \, \treenv{\vertex10} + 3 \ \treenv{}^3.
\end{align*}
In these examples, the coefficient $2$ (resp.~$3$) before a forest indicates that there are two (resp.~three) mappings that yield the forest. 
In addition, $ \treenv{}^2 $ denotes the forest $ \treenv{}\, \treenv{} $.
Note that the pruning $ \tree \setminus \tree' $ is an element of the forest space $\ForestS$.

\section{Proposition~\ref{prop1} and the proof of Lemma~\ref{lemma1}}

The following proposition is a special case $\tree' = \tau$ in Lemma~\ref{lemma1}.
We show its proof for the readers' convenience as it is not 
given in~\cite{bu21}.

\begin{prop}[\protect{cf.~\cite[Corollary~3.4F]{bu21}}]\label{prop1}
    Let $a\in\Bg$. 
    Then, for any positive integer $n$, it holds that
    \[ f^{(n)} \left( (\BB_h a ) y_0 \right) = \sum_{ \tree \in \Tree } h^{ | \tree | - 1 } \sigma (\tree)^{-1} (a \D) (\tree) \F ( \tree \ast^n ) = \sum_{ \forest \in \Forest } h^{ | \forest | } \sigma ([\forest])^{-1} (a \D) ([\forest]) \F ( [\forest] \ast^n ). \]
\end{prop}

\begin{proof}
    By introducing $ \boldK = \{ \boldk = [ k_1, k_2 , \dots, k_m] \in \mathbb{Z}_{\ge0}^m \mid m \in \mathbb{Z}_{> 0} \} $, Taylor's theorem implies that
    \[ f^{(n)} \bigg( y_0 + \sum_{i=1}^{\infty} \delta_i \bigg) = \sum_{ \boldk \in \boldK } \frac{1}{\boldk !} \f^{ ( n + | \boldk | )} \delta_1^{k_1} \cdots \delta_m^{k_m}, \]
    where $ \boldk ! \coloneqq \prod_{i=1}^m k_i ! $ and $ |\boldk| = \sum_{i=1}^m k_i $.
    
    By setting $ \delta = (\BB_h a ) y_0 - y_0 $ in the above equation, we have
    \begin{align*}
        &f^{(n)} \left( (\BB_h a ) y_0 \right)\\
        &\quad = \sum_{ \boldk } \frac{1}{\boldk !} \sum_{ \tree_1 \in \Tree } \cdots \sum_{ \tree_m \in \Tree } \frac{ h^{ \sum k_i | \tree_i| } }{ \sigma (\tree_1)^{k_1} \cdots \sigma (\tree_m)^{k_m} } a (\tree_1)^{k_1} \cdots a (\tree_m)^{k_m} \f^{ (n+|\boldk|) } \F (\tree_1)^{k_1} \cdots \F (\tree_m)^{k_m} \\
        &\quad = \sum_{\tree \in \Tree } \frac{ h^{|\tree|-1} }{\sigma (\tree)} a \D (\tree) \F (\tree \ast^n ).
    \end{align*}
    This establishes the first equality in the proposition. 
    The second equality is straightforward since $ [\cdot ] \colon \Forest \to \Tree $ is bijective. 
\end{proof}

We now give a new proof for Lemma~\ref{lemma1}.

\begin{proof}

The case $\tree' = \tau$ is covered by Proposition~\ref{prop1}. %, whose proof is also given in the next section.
The general case will be proved by induction with respect to the operator $ [\cdot ] $.
Below, let $ y_1 = \left( \BB_h a \right) y_0 $.

% We prove the lemma using the induction with respect to the operator $ [\cdot ] $. 
% For $\tree = \tau$, the claim (Proposition~\ref{prop1}) is easy to verify.
% The proof is given in the next section.

Consider the tree $ \tree' = \left[ \prod^{n'}_{i=1} \left( \tree'_i \right)^{k'_i} \right]$, where $ \tree'_1, \tree'_2, \dots,\tree'_{n'} $ are distinct.
Let $ N' = \sum_{i=1}^{n'} k'_i $.
Using Proposition~\ref{prop1} and the induction hypothesis, we obtain
\begin{align*}
    &\frac{ h^{|\tree'|} }{\sigma (\tree')} \F (\tree') \left( y_1 \right) \\
    &\quad = \frac{ h^{|\tree'|} }{\sigma (\tree')} f^{ (N') } ( y_1 ) \prod_{i=1}^{n'} \left( \F \left( \tree'_i \right) \left( y_1 \right) \right)^{k'_i} \\
    &\quad = \frac{ h^{|\tree'|} }{\sigma (\tree')} \sum_{\forest \in \Forest} \frac{ h^{|\forest|} }{\sigma \left( [\forest] \right)} a \D \left( [\forest] \right) \F \left( [\forest] \ast^{N'} \right) \prod_{i=1}^{n'} \left( \F \left( \tree'_i \right) \left( y_1 \right) \right)^{k'_i} \\
    &\quad = \frac{ h^{|\tree'|} }{\sigma (\tree')} \sum_{\forest \in \Forest} \frac{ h^{|\forest|} }{\sigma \left( [\forest] \right)} a \D \left( [\forest] \right) \F \left( [\forest] \ast^{N'} \right) \prod_{i=1}^{n'} \prod_{j=1}^{k'_i} \sum_{\tree_{ij} \ge \tree'_i} \frac{ h^{|\tree_{ij}|} \sigma \left( \tree'_i \right) }{ h^{|\tree'_i|} \sigma \left(\tree_{ij} \right) } a \left( \tree_{ij} \setminus \tree'_i \right) \F \left( \tree_{ij} \right).
\end{align*}
Now, define the set
\begin{equation*}
    S (\tree') = \left\{ \left(\forest, \{ \tree_{ij}\} \right) : \forest \in \Forest, \ \tree_{ij} \ge \tree'_i \ (i=1,2,\dots,n'; \, j = 1,2,\dots, k'_i) \right\}
\end{equation*}
and let $ \tree = \left[ \forest \prod_{i,j} \tree_{ij} \right] $.
Then, we have 
\begin{equation}\label{eq:LHSthm39B}
    \frac{ h^{|\tree'|} }{\sigma (\tree')} \F (\tree') \left( y_1 \right)
    = \sum_{ \left(\forest, \{ \tree_{ij}\} \right) \in S (\tree') } \frac{ h^{|\tree|} }{ \left(\prod_{i=1}^{n'} k'_i ! \right) \sigma \left( [\forest] \right) \prod_{i,j} \sigma \left( \tree_{ij} \right) } a ( \forest ) \prod_{i,j} a \left( \tree_{ij} \setminus \tree'_i \right) \F \left( \tree \right).
\end{equation}
Here, we have used the relation $\sigma(\tree') = \prod_{i=1}^{n'}(k_i'!) (\sigma(k_i'))^{k_i'}$.

Next, we consider the right-hand side of Lemma~\ref{lemma1}.
To deal with the pruning $ \tree \setminus \tree' $, where $ \tree = \left[ \prod^{n}_{j=1} \left( \tree_j \right)^{k_j} \right]$ and $ \tree' = \left[ \prod^{n'}_{i=1} \left( \tree'_i \right)^{k'_i} \right]$
with distinct trees $ \tree_1, \tree_2, \dots,\tree_{n} $ and $ \tree'_1, \tree'_2, \dots,\tree'_{n'} $ (note that $\tree_i$ can be equal to $\tree'_j$),
we introduce the following definition. 
An $ (n'+1) \times n$ non-negative integer matrix $ m $ (the $(i,j)$ entry of $m$ is denoted by $m_{ij}$, where $ i = 0,1,...,n'$ and $j = 1,2,...,n$)
is said to be an assignment for the pair $ (\tree',\tree) $ if it satisfies the following three conditions:
\begin{itemize}[(iii)]
    \item[(i)] if $ m_{ij} > 0 $, then $ \tree_j \ge \tree'_i $, 
    \item[(ii)] for all $ j $, $ \sum_{i=0}^{n’} m_{ij} = k_j $, and
    \item[(iii)] for all $ i \neq 0$, $ \sum_{j=1}^{n} m_{ij} = k’_i $.
\end{itemize}
Each assignment $ m $ corresponds to a specific pruning operation given by
\[ \prod_{j=1}^n \left( \tree_j^{ m_{0j} } \prod_{i=1}^{n'} \left( \tree_j \setminus \tree'_i \right)^{m_{ij}} \right) = \prod_{j=1}^n \prod_{i=0}^{n'} \left( \tree_j \setminus \tree'_i \right)^{m_{ij}}, \]
where $ t'_0 = \emptytree $ and we have adopted the convention $ \tree \setminus \emptytree = \tree $.
Let $ M (\tree',\tree) $ be the set of all assignments for $ (\tree',\tree) $. 
Then, we have
\begin{equation*}
    \tree \setminus \tree'
    = \sum_{ m \in M (\tree', \tree) } \prod_{j=1}^n \left( \frac{ k_j ! }{ \prod_{i=0}^{n'} m_{ij} ! } \prod_{i=0}^{n'} \left( \tree_j \setminus \tree'_i \right)^{m_{ij}} \right).
\end{equation*}
Here, the multinomial coefficient represents the number of combinations
to assign the trees $ \tree_j^{k_j} $ to the trees $ \{ \left(\tree'_i\right)^{m_{ij}} \}_{i=0}^{n'} $ 
(see Remark~\ref{rem:ex-assignment} for an illustrating example of the relation between assignments and pruning).
Using this expression of pruning and letting $ \forest = \prod_{j=1}^n \tree_j^{ m_{0j} } $, the summed expression of the right-hand side of Lemma~\ref{lemma1} can be written as
\begin{align*}
    &\frac{h^{|\tree|}}{ \sigma \left( \tree \right) } a \left( \tree \setminus \tree' \right) \F (\tree) \\
    &\quad = \frac{h^{|\tree|}}{ \prod_{j=1}^n k_j ! \left( \sigma \left( \tree_j \right) \right)^{k_j} } \\
    &\qquad \cdot \sum_{ m \in M (\tree', \tree) } \prod_{j=1}^n \left( \frac{ k_j ! }{ \prod_{i=0}^{n'} m_{ij} ! } a \left( \tree_j^{ m_{0j} } \right) \prod_{i=1}^{n'} a \left( \tree_j \setminus \tree'_i \right)^{m_{ij}} \right) \cdot \F (\tree) \\
    &\quad = \sum_{ m \in M (\tree', \tree) } \frac{h^{|\tree|}}{\sigma ([\forest]) \prod_{i=1}^{n'} \prod_{j=1}^n \left( m_{ij} ! \left( \sigma \left( \tree_j \right) \right)^{m_{ij}} \right) } \cdot a (\forest) \prod_{i=1}^{n'} \prod_{j=1}^n a \left(\tree_j \setminus \tree'_i \right)^{m_{ij}} \cdot \F (\tree).
\end{align*}
On the other hand, in the summation in \eqref{eq:LHSthm39B}, the term which corresponds to an assignment $ m $ appears $ \prod_{i=1}^{n'} \frac{k'_i !}{ \prod_{j=1}^n m_{ij} ! } $ times. 
Therefore, we have 
\begin{align*}
    &\frac{ h^{|\tree'|} }{\sigma (\tree')} \F (\tree') \left( y_1 \right) \\
    &\quad = \sum_{ \left(\forest, \{ \tree_{ij}\} \right) \in S (\tree') } \frac{ h^{|\tree|} }{ \left(\prod_{i=1}^{n'} k'_i ! \right) \sigma \left( [\forest] \right) \prod_{i=1}^{n'} \prod_{j=1}^n \sigma \left( \tree_{ij} \right) } a ( \forest ) \prod_{i=1}^{n'} \prod_{j=1}^n a \left( \tree_{ij} \setminus \tree'_i \right) \F \left( \tree \right)\\
    &\quad = \sum_{ \tree \ge \tree'} \sum_{ m \in M (\tree',\tree) }  \frac{ \prod_{i=1}^{n'} k'_i !}{ \prod_{i=1}^{n'} \prod_{j=1}^n m_{ij} ! } \frac{ h^{|\tree|} }{ \left(\prod_{i=1}^{n'} k'_i ! \right) \sigma \left( [\forest] \right) \prod_{i=1}^{n'} \prod_{j=1}^n \left( \sigma \left( \tree_j \right) \right)^{m_{ij}} } \\
    &\qquad \qquad \cdot a ( \forest ) \prod_{i=1}^{n'} \prod_{j=1}^n a \left(\tree_j \setminus \tree'_i \right)^{m_{ij}} \cdot \F \left( \tree \right)\\
    &\quad = \sum_{ \tree \ge \tree'} \frac{h^{|\tree|}}{ \sigma \left( \tree \right) } a \left( \tree \setminus \tree' \right) \F (\tree),
\end{align*}
which proves the lemma.
\end{proof}

\begin{rem}\label{rem:ex-assignment}
    In this remark, we explain how the concept of assignment can be used to consider pruning, using the following example: 
    \begin{align*}
        \tree &= \treenv{\vertex1{-1.5}\vertex2{-0.75}\vertex10\vertex40\vertex10\vertex6{0.75}\vertex10\vertex8{1.8}\vertex1{-0.3}\vertex2{0.3}},&
        \tree' &= \treenv{\vertex1{-0.8}\vertex20\vertex10\vertex4{0.8}\vertex10}.
    \end{align*}
    In this case, these trees can be expressed as $ \tree = [ \tree_1^{k_1} \tree_2^{k_2}\tree_3^{k_3} ] $ and $ \tree' = [ (\tree'_1)^{k'_1} (\tree'_2)^{k'_2}] $, where 
    \begin{align*}
        \tree_1 &= \treenv{},&
        \tree_2 &= \treenv{\vertex10},&
        \tree_3 &= \treenv{\vertex1{-0.6}\vertex2{0.6}},&\
        \tree'_1 &= \treenv{},&
        \tree'_2 &= \treenv{\vertex10},\\
        k_1 &= 1,&
        k_2 &= 3,&
        k_3 &= 1,&
        k'_1 &= 1,&
        k'_2 &= 2.
    \end{align*}
    Since $ \tree_1 \ge \tree'_2 $ does not hold, the condition (i) of assignments implies that $ m_{21} = 0 $. 
    Consequently, it is easy to verify
    \[ M (\tree',\tree) = \left\{
    \begin{bmatrix} 0 & 1 & 1 \\ 1 & 0 & 0 \\ 0 & 2 & 0 \end{bmatrix},
    \begin{bmatrix} 0 & 2 & 0 \\ 1 & 0 & 0 \\ 0 & 1 & 1 \end{bmatrix},
    \begin{bmatrix} 1 & 0 & 1 \\ 0 & 1 & 0 \\ 0 & 2 & 0 \end{bmatrix},
    \begin{bmatrix} 1 & 1 & 0 \\ 0 & 1 & 0 \\ 0 & 1 & 1 \end{bmatrix},
    \begin{bmatrix} 1 & 1 & 0 \\ 0 & 0 & 1 \\ 0 & 2 & 0 \end{bmatrix}
    \right\}. \]
    
    The assignment
    \begin{equation*}
        m =
        \begin{bNiceMatrix}[first-row,first-col]
            & \treenv{} & \treenv{\vertex10} & \treenv{\vertex1{-0.6}\vertex2{0.6}} \\
            \emptytree & 0 & 1 & 1 \\
            \raisebox{3pt}{\treenv{\frame{-0.4,-0.5}{0.6,1.2}}} & 1 & 0 & 0 \\
            \treenv{\frame{-0.4,-0.5}{0.6,1.2}\vertex10} & 0 & 2 & 0
        \end{bNiceMatrix}
    \end{equation*}
    corresponds to the case that
    ``$ \tree'_1 $ is regarded as a subtree of $ \tree_1 $, both of two copies of $ \tree'_2 $ are regarded as subtrees of two copies of $ \tree_2 $, and $ \tree_3 $ and one copy of $ \tree_2 $ remain untouched''; in other words, remaining vertices and edges are expressed as $ \treenv{\vertex10} \ \treenv{\vertex1{-0.6}\vertex2{0.6}} $. 
    Here, since there are three copies of $ \tree_2 $, there are three possible ways to choose two of them. 
    As a consequence, this assignment corresponds to the term $ 3 \, \treenv{\vertex10} \ \treenv{\vertex1{-0.6}\vertex2{0.6}} $, which is appeared in the pruning $ \tree \setminus \tree' $. 

    Considering the other assignments in a similar way, we obtain the following results: 
    \begin{align*}
        \begin{bmatrix} 0 & 2 & 0 \\ 1 & 0 & 0 \\ 0 & 1 & 1 \end{bmatrix}& \rightarrow \ 3 \, \treenv{\vertex10}^2 \left(\treenv{\vertex1{-0.6}\vertex2{0.6}} \setminus \treenv{\vertex10}\right) = 6 \, \treenv{} \, \treenv{\vertex10}^2 ,&
        \begin{bmatrix} 1 & 0 & 1 \\ 0 & 1 & 0 \\ 0 & 2 & 0 \end{bmatrix}& \rightarrow \ 3 \, \treenv{} \left( \treenv{\vertex10} \setminus \treenv{} \right) \treenv{\vertex1{-0.6}\vertex2{0.6}} = 3 \, \treenv{}^2 \, \treenv{\vertex1{-0.6}\vertex2{0.6}} ,\\
        \begin{bmatrix} 1 & 1 & 0 \\ 0 & 1 & 0 \\ 0 & 1 & 1 \end{bmatrix}& \rightarrow \ 6 \, \treenv{} \, \treenv{\vertex10} \left( \treenv{\vertex10} \setminus \treenv{} \right)  \left(\treenv{\vertex1{-0.6}\vertex2{0.6}} \setminus \treenv{\vertex10}\right) = 12 \, \treenv{}^3 \,\treenv{\vertex10} ,&
        \begin{bmatrix} 1 & 1 & 0 \\ 0 & 0 & 1 \\ 0 & 2 & 0 \end{bmatrix}& \rightarrow \ 3 \, \treenv{} \, \treenv{\vertex10} \left(\treenv{\vertex1{-0.6}\vertex2{0.6}} \setminus \treenv{}\right) = 3 \, \treenv{}^3 \, \treenv{\vertex10} .
    \end{align*}
    Adding these expressions, we obtain 
    \[ \tree \setminus \tree' = 15 \, \treenv{}^3 \, \treenv{\vertex10} + 3 \, \treenv{}^2 \, \treenv{\vertex1{-0.6}\vertex2{0.6}} + 6 \, \treenv{} \, \treenv{\vertex10}^2 + 3 \, \treenv{\vertex10} \ \treenv{\vertex1{-0.6}\vertex2{0.6}}. \]
\end{rem}

\backmatter

% \bmhead{Supplementary information}

% If your article has accompanying supplementary file/s please state so here. 

% Authors reporting data from electrophoretic gels and blots should supply the full unprocessed scans for key as part of their Supplementary information. This may be requested by the editorial team/s if it is missing.

% Please refer to Journal-level guidance for any specific requirements.

% \bmhead{Acknowledgements}

% Acknowledgements are not compulsory. Where included they should be brief. Grant or contribution numbers may be acknowledged.

% Please refer to Journal-level guidance for any specific requirements.

\section*{Declarations}

% Some journals require declarations to be submitted in a standardised format. Please check the Instructions for Authors of the journal to which you are submitting to see if you need to complete this section. If yes, your manuscript must contain the following sections under the heading `Declarations':

\begin{itemize}
\item Funding\\
The third author was supported by JSPS Grant-in-Aid for Scientific Research (A) (No.~20H00581), JSPS Grant-in-Aid for Challenging Research (Pioneering) (No.~21K18301), JSPS 	
Grant-in-Aid for Scientific Research (B) (No.~20H01822, No.~24K02951 and No.~24K00540) and JST PRESTP Grand No. JPMJPR2129.
The fourth author was supported by JSPS Grant-in-Aid for Early-Career Scientists (No.~22K13955), JSPS Grant-in-Aid for Scientific Research (B) (No.~20H01822 and No.~24K00540), JSPS Grant-in-Aid for Challenging Research (Exploratory) (No.~24K22290). 
\item Competing interests\\
Not applicable
\item Ethics approval and consent to participate\\
The key idea of the new proof was already published in Japanese by Maruzen Publisher in~\cite{bu24} (note that it was not included in~\cite{bu21}). 
In this paper, we have reorganized the proof, incorporating further corrections and new examples not present in~\cite{bu24}.
We have obtained permission from the publisher to present this work in English and submit it to a scientific journal.
% \item Consent for publication\\
% Not applicable
% \item Data availability\\
% Not applicable
% \item Materials availability\\
% Not applicable
% \item Code availability\\
% Not applicable
\end{itemize}

\end{document}